\documentclass[12pt]{article}
\usepackage{amsmath}
\usepackage{amssymb}
\usepackage{amscd}
\def\G{{\rm G}} %the grassmanian
\def\C{{\cal C}}  %the class C
\newtheorem{lemme}{Lemme}
\newtheorem{theoreme}{Th\'eor\`eme}
\newtheorem{proposition}{Proposition}
\newtheorem{corollaire}{Corollaire}
\newcounter{numeroexemple}
\newenvironment{exemple}
  {\addtocounter{numeroexemple}{1} 
  \begin{trivlist}\item[]\textbf{Exemple \thenumeroexemple}}{\end{trivlist}}
\newenvironment{remarque}
  {\begin{trivlist}\item[]\textbf{Remarque}}{\end{trivlist}}
\newcounter{numerodefinition}

\newcounter{numeroquestion}
\newenvironment{question}
  {\addtocounter{numeroquestion}{1}
\begin{trivlist}\item[]\textit{Question \thenumeroquestion}}
{\end{trivlist}}
\newcounter{numeroconjecture}
\newenvironment{conjecture}
{  \addtocounter{numeroconjecture}{1}
\begin{trivlist}\item[]\textit{Conjecture \thenumeroconjecture}}
{\end{trivlist}}
\newenvironment{preuve}{\begin{trivlist}\item[]\textit{Preuve.}}
{\item[] \hfill $\square$\end{trivlist}}
\newcommand{\Tan}{{\rm Tan\,}}
\begin{document}
\title{Sur la caract\'erisation du bord d'une cha\^{\i}ne holomorphe
  dans l'espace projectif}
\author{{\sc Tien-Cuong} DINH\footnote{
Math\'ematique-B\^atiment 425,
Universit\'e Paris-Sud, 91405 ORSAY Cedex (France).\break
Mots cl\'es: probl\`eme du bord, 
conditions des moments, maximalement
complexe.\break 
Classification AMS: 32C25, 32C16.
}}
\date{}
\maketitle
{\small
\noindent
{\bf R\'esum\'e:}  Nous d\'emontrons qu'une
  sous-vari\'et\'e r\'eelle, compacte, orient\'ee et lisse $\Gamma$ de
  dimension $2p-1\geq 3$ de $\mathbb{CP}^n$ est le bord d'un
  sous-ensemble analytique
  s'il existe une vari\'et\'e r\'eelle 
$V\subset\G(n-p+2,n+1)$ de codimension $1$
 satisfaisant les conditions suivantes pour tout
  $\nu\in V$
\begin{enumerate}
\item La r\'eunion $\bigcup_{\nu\in V}\mathbb{P}^{n-p+1}_\nu$ recouvre
  un ouvert dense de $\Gamma$. 
\item Le $(n-p+1)$-plan 
$\mathbb{P}^{n-p+1}_\nu$ intersecte $\Gamma$ transversalement.
\item $\Gamma\cap\mathbb{P}^{n-p+1}_\nu$ est le bord d'une surface de
  Riemann dans $\mathbb{P}^{n-p+1}_\nu$.
\item Aucun ouvert non vide de 
$\Gamma\cap\mathbb{P}^{n-p+1}_\nu$ n'est r\'eel analytique.
\end{enumerate}
Pour la preuve, nous utilisons et d\'emontrons le r\'esultat suivant:  
pour toute surface de Riemann \`a bord rectifiable
  (\'eventuellement r\'eductible et singuli\`ere) $S$ d'une vari\'et\'e
  complexe, $\overline S$ admet un syst\`eme fondamental 
de voisinages de Stein. Il existe $\Gamma$ r\'eelle alg\'ebrique
v\'erifiant 1-3, qui n'est pas bord d'un sous-ensemble analytique.
}
\\
\section{R\'esultats} 
L'origine du probl\`eme de caract\'erisation du bord d'une vari\'et\'e
complexe
ou d'une cha\^{\i}ne holomorphe (probl\`eme du bord) 
est un th\'eor\`eme de Wermer qui
affirme qu'une courbe r\'eelle, analytique, ferm\'ee et orient\'ee 
$\gamma$ dans
$\mathbb{C}^n$ est le bord d'une surface de Riemann $S$ si et seulement si
elle v\'erifie {\it la condition des moments},
c.-\`a.-d. $\int_\gamma\varphi=0$ pour toute $(1,0)$-forme $\varphi$
holomorphe dans $\mathbb{C}^n$. Cette condition est r\'eduite de la
formule de Stokes sur $\overline S$. Dans le cas o\`u la condition des
moments est valide, l'enveloppe polynomiale
de $\gamma$ est \'egale \`a la r\'eunion $S\cup\gamma$. Dans le cas
contraire, l'enveloppe polynomiale
de $\gamma$ est \'egale \`a elle-m\^eme. L'enveloppe polynomiale d'un
compact de longueur finie a \'et\'e \'etudi\'ee ensuite dans plusieurs
travaux de Bishop \cite{Bishop}, de Stolzenberg \cite{Stolzenberg},
d'Alexander \cite{Alexander1}, Lawrence \cite{Lawrence1, Lawrence2} et
\'egalement dans \cite{Dinh2, Dinh3}. Les r\'esultats de ces travaux 
sont utilisables pour
d\'emontrer de nouveau 
le th\'eor\`eme de Harvey-Lawson, qui donne la solution
du probl\`eme du bord en dimension sup\'erieure \cite{HarveyLawson}. 
Ce th\'eor\`eme de
Harvey-Lawson dit qu'une sous-vari\'et\'e r\'eelle, compacte,
orient\'ee $\Gamma$ de
dimension $2p-1\geq 3$ de $\mathbb{C}^n$ est le bord d'une vari\'et\'e
complexe si et seulement si elle est {\it maximalement complexe},
c.-\`a-d. son plan tangent en chaque point contient un sous-espace
complexe de dimension $p-1$ (la dimension maximale possible). La
condition ``$\Gamma$ maximalement complexe'' implique que la courbe
$\Gamma\cap\mathbb{C}^{n-p+1}_\nu$ v\'erifie la condition des moments
pour toute $(n-p+1)$-plan $\mathbb{C}^{n-p+1}_\nu$ qui coupe $\Gamma$
transversalement. Dans $\mathbb{CP}^n$ le th\'eor\`eme de
Harvey-Lawson n'est plus valable. D'apr\`es le th\'eor\`eme de
Dolbeault-Henkin g\'en\'eralis\'e, si $\Gamma$ est maximalement
complexe et si l'intersection $\Gamma\cap\mathbb{P}^{n-p+1}_\nu$ borde
une surface de Riemann pour une famille 
suffisamment grande $V$ de $(n-p+1)$-plans projectifs, alors $\Gamma$ borde
une vari\'et\'e complexe \cite{DolbeaultHenkin, Dinh4}. En
particulier,  ce th\'eor\`eme est valable pour 
une vari\'et\'e r\'eelle $V$ de codimension $1$ de la
grassmannienne $\G(n-p+2,n+1)$. Cette grassmannienne 
est l'espace de param\`etre  des
$(n-p+1)$-plans projectifs dans $\mathbb{CP}^n$. 
Nous allons d\'emontrer que dans le cas
pr\'ec\'edent, sous une hypoth\`ese sur $\Gamma$ (qui est valable
g\'en\'eriquement), si $\bigcup_{\nu\in V} \mathbb{P}^{n-p+1}_\nu$
recouvre $\Gamma$, l'hypoth\`ese ``$\Gamma$ maximalement complexe'' ne
sera plus n\'ecessaire. Le probl\`eme du bord sans l'hypoth\`ese
``$\Gamma$ maximalement complexe'' 
a \'et\'e pos\'e dans $\mathbb{C}^n$ par
Dolbeault-Henkin et \'egalement 
par Globevnik-Stout sous la forme d'un probl\`eme
d'extension holomorphe \cite{GlobevnikStout}, 
dont les solutions pour un cas g\'en\'erique
se trouvent dans \cite{Dinh5}. La preuve des r\'esultats cit\'es
ci-dessus utilise essentiellement la condition des moments et elle
est extensible au cas dans $\mathbb{CP}^n$ gr\^ace au th\'eor\`eme de
Mihalache g\'en\'eralis\'e (th\'eor\`eme 2). \\
On appelle {\it ouvert $(n-p+1)$-lin\'eairement concave } tout ouvert
$X$ de $\mathbb{CP}^n$ qui est la r\'eunion d'une famille continue de
$(n-p+1)$-plans projectifs. On pose $X^*:=\{\nu\in \G(n-p+2,n+1)\mbox{
  telle que }\mathbb{P}^{n-p+1}_\nu\subset X\}$. C'est un ouvert
connexe. Pour tout ferm\'e $\Gamma$ de $X$, {\it une $p$-cha\^{\i}ne
  holomorphe} de $X\setminus\Gamma$ est une combinaison lin\'eaire,
localement finie \`a coefficients entiers de sous-ensembles
analytiques de dimension pure $p$ de $X\setminus\Gamma$. Toute
$p$-cha\^{\i}ne holomorphe de $X\setminus\Gamma$ est de mesure
$2p$-dimensionnelle 
localement finie dans $X$ \cite{Dinh4}.
\begin{theoreme} Soient $X$ un ouvert $(n-p+1)$-lin\'eairement concave
  de $\mathbb{CP}^n$ et $\Gamma$ une combinaison lin\'eaire, finie,
  disjointe \`a
  coefficients entiers de sous-vari\'et\'es r\'eelles, $C^2$ de
  dimension $2p-1\geq 3$ de $X$. Soit $V$ une vari\'et\'e
  r\'eelle de codimenion $1$ de $X^*$ v\'erifiant les conditions
  suivantes
\begin{enumerate}
\item  La r\'eunion $\bigcup_{\nu\in V}\mathbb{P}^{n-p+1}_\nu$ recouvre un
  ouvert dense de
  $\Gamma$.
\item L'intersection $\Gamma\cap \mathbb{P}^{n-p+1}_\nu$ est
  transversale pour tout $\nu\in V$.
\item L'intersection $\Gamma\cap \mathbb{P}^{n-p+1}_\nu$ est le bord
  d'une $1$-cha\^{\i}ne holomorphe de
  $\mathbb{P}^{n-p+1}_\nu\setminus\Gamma$ au sens des courants.
\item Aucun ouvert non vide de $\Gamma\cap
  \mathbb{P}^{n-p+1}_\nu$ n'est r\'eel analytique.
\end{enumerate}
Alors $\Gamma$ est le bord d'une $p$-cha\^{\i}ne holomorphe de
$X\setminus \Gamma$, de masse localement finie dans $X$ au sens des
courants.
\end{theoreme}
\begin{remarque} Ce th\'eor\`eme reste valable pour une famille plus large
des ensembles $V$, qui ne sont pas n\'ecessairement des
  vari\'et\'es. Pour un r\'esultat analogue dans $X=\mathbb{C}^n$, 
cette famille est d\'efinie pr\'ecis\'ement dans \cite{Dinh5}.
\end{remarque}
L'exemple suivant montre que la condition 4 dans le th\'eor\`eme
pr\'ec\'edent est n\'ecessaire:
\begin{exemple} (Henkin \cite{Dinh5}) 
Soit $\Gamma\subset\mathbb{C}^{3}\subset\mathbb{CP}^{3}$
une vari\'et\'e r\'eelle alg\'ebrique de dimension $3$ d\'efinie par
$$
\Gamma=
\{ y_{2}=y_{3}=0,x_{1}^{2}+y_{1}^{2}+x_{2}^{2}+x_{3}^{2}=1
\}
$$
o\`u $z_{1}=x_{1}+iy_{1}$, $z_{2}=x_{2}+iy_{2}$, $z_{3}=x_{3}+iy_{3}$
sont les coordonn\'ees de $\mathbb{C}^{3}$. Consid\'erons l'hyperplan
projectif 
$$H_{a,b,c}=\{z_{1}=az_{2}+bz_{3}+c\}$$ 
o\`u $a=a_{1}+ia_{2}$,
 $b=b_{1}+ib_{2}$ et $c=c_{1}+ic_{2}$. Posons $\Gamma_{a,b,c}=\Gamma\cap
H_{a,b,c}$. Alors pour un $(a,b,c)$ g\'en\'erique,
$\Gamma_{a,b,c}$ est une courbe r\'eelle ferm\'ee de la surface de Riemann
alg\'ebrique $S_{a,b,c}\subset H_{a,b,c}$ qui est d\'efinie par
$$S_{a,b,c}=\{(a_{1}z_{2}+b_{1}z_{3}+c_{1})^{2} +
(a_{2}z_{2}+b_{2}z_{3}+c_{2})^{2} +z_{2}^{2}+z_{3}^{2}=1\}\cap H_{a,b,c}.
$$
Comme $S_{a,b,c}$ est une surface de Riemann compacte 
de genre $0$, la courbe $\Gamma\cap H_{a,b,c}$ borde
une surface de Riemann dans $\mathbb{CP}^{3}$. 
La vari\'et\'e $\Gamma$ ne peut pas \^etre le
bord d'une vari\'et\'e complexe car elle n'est pas maximalement
complexe.
\end{exemple}
Soit $S$ un sous-espace analytique d'une vari\'et\'e complexe $Z$. Si
$S$ est de Stein, d'apr\`es le th\'eor\`eme de Siu, $S$ admet un
syst\`eme fondamental de voisinages de Stein. En particulier, si $S$
est une surface de Riemann ouverte, dont le bord est une r\'eunion finie de
courbes de Jordan de longueur finie, alors $S$ admet
dans $Z$ un syst\`eme fondamental de voisinages de Stein. On appelle {\it 
  vari\'et\'e de Stein} toute vari\'et\'e plongeable dans une
vari\'et\'e affine complexe. Un {\it
  voisinage de Stein} de $S$ est un ouvert contenant $S$ tel que
chacune de ses composantes connexes soit de Stein. 
Le
th\'eor\`eme suivant donne la r\'eponse \`a une question de
Dolbeault-Henkin \cite{Mihalache}:
\begin{theoreme} Soient $Z$ une vari\'et\'e complexe, $\Gamma\subset Z$
  une r\'eunion finie, disjointe 
de courbes de Jordan, de longueur finie et $S$ un sous
ensemble analytique de dimension pure $1$ de 
$Z\setminus \Gamma$, born\'e dans $Z$. Supposons que $S\cup\Gamma$ ne
contient aucun sous-ensemble analytique, compact, de dimension pure
$1$ de $Z$. Alors 
\begin{enumerate}
\item $S\cup\Gamma$ admet un voisinage connexe de Stein.
\item $S\cup\Gamma$ admet un syst\`eme fondamental de voisinages de Stein.
\end{enumerate}
\end{theoreme}
\begin{remarque}
 Si $Z=\mathbb{C}^n$ et si $\Gamma$ est une courbe  irr\'eductible,
 alors $\overline S$ est \'egal \`a {\it l'enveloppe polynomiale} 
$$\widehat\Gamma
:=\left\{x\in\mathbb{C}^n:\ |P(x)|\leq \max_{z\in\Gamma}
 |P(z)| \mbox{ pour tout polyn\^ome } P\right\}$$ de $\Gamma$ 
\cite{Wermer, Bishop, Stolzenberg, Alexander1}. D'apr\`es un
 lemme d'Oka, $\overline S$ admet un syst\`eme fondamental de
 voisinages de Stein de la forme $\{|P_k(z)|<1;\ k=1,\ldots,m\}$, o\`u
 $P_k$ sont des polyn\^omes. En g\'en\'eral, dans $\mathbb{C}^n$,
 Henkin \cite{Fabre} et Alexander-Wermer \cite{AlexanderWermer} 
ont prouv\'e ind\'ependamment que
 $\overline S$ est rationnellement convexe, c.-\`a-d. que
 $\mathbb{C}^n\setminus S$ est une r\'eunion d'hypersurfaces
 alg\'ebriques de $\mathbb{C}^n$. Par
 cons\'equent, $\overline S$ admet un syst\`eme de voisinages de Stein
 de la forme $\{|R_k|<1 \mbox{ pour } k=1,\ldots,m\}$, o\`u $R_k$ sont
 des fonctions rationnelles, holomorphes au voisinage de $\overline
 S$.\\
Dans $Z=\mathbb{CP}^n$, Fabre a construit une surface de Riemann
 irr\'eductible, \`a
 bord lisse, non rationnellement convexe \cite{Fabre}. Si $\Gamma$ est
 irr\'eductible dans $Z=\mathbb{CP}^n$, 
le th\'eor\`eme 2 se r\'eduit \`a celui de Mihalache
 \cite{Mihalache}.
\end{remarque}
 La d\'emonstration de Mihalache est une application
 du th\'eor\`eme de Siu et celui de Stolzenberg-Alexander
 \cite{Alexander1}. Notre d\'emonstration est un prolongement du
 travail de Mihalache. Nous utilisons  le
 th\'eor\`eme de Henkin-Alexander-Wermer pour compl\'eter la preuve et
 appliquer ce r\'esultat pour d\'emontrer les th\'eor\`emes
  1 et 3.\\
\section{Preuve du th\'eor\`eme 2}
\begin{lemme}\textbf{(\cite{Alexander1})} 
Soient $Z$ une vari\'ete
 de Stein, $K\subset Z$ un compact holomorphiquement convexe et
 $\Gamma\subset Z$ une r\'eunion finie d'arcs r\'eels compacts de
 longueur finie. Alors $(K\cup\Gamma)^\wedge\setminus(K\cup\Gamma)$
 est un sous-ensemble analytique de dimension pure 1 (\'eventuellement
 vide) de $Z\setminus (K\cup\Gamma)$, born\'e dans $Z$, o\`u
 $$(K\cup\Gamma)^\wedge:= \{z\in Z:\ |h(z)|\leq
 \max_{x\in K\cup\Gamma}|h(x)| \mbox{ pour } h \mbox{
 holomorphe dans } Z\}$$
est l'enveloppe d'holomorphie de $K\cup\Gamma$ dans $Z$.
\end{lemme}
\begin{lemme}\textbf{(\cite{Fabre,AlexanderWermer})} Sous l'hypoth\`ese du
 th\'eor\`eme 2, si $Z$ est de Stein, $S\cup\Gamma$ sera
 m\'eromorphiquement convexe, c.-\`a.-d. $Z\setminus(S\cup\Gamma)$ est
 une r\'eunion d'hypersurfaces complexes de $Z$.
\end{lemme}
\begin{preuve} Si $Z$ est de Stein, elle est plongeable dans un
 $\mathbb{C}^N$. La preuve du lemme pr\'ec\'edent se ram\`ene au cas
 dans $\mathbb{C}^N$, qui se trouve dans \cite[lemme 1.5]{AlexanderWermer}. 
\end{preuve}
\begin{corollaire} Soit $S_i$ une suite de surfaces de Riemann
 born\'ees 
dans $\mathbb{C}^n$. Supposons que pour $i=0,1,\ldots$ le bord de $S_i$
  est une r\'eunion finie, disjointe  de courbes de Jordan
  rectifiables. Supposons de plus que les
  bords des $S_i$ tendent vers le bord de $S_0$ au sens g\'eom\'etrique et
  aussi au sens des courants. Alors $\overline S_i$ tend vers
  $\overline S_0$ au sens g\'eom\'etrique. 
\end{corollaire}
\begin{preuve}
Soit $z\not \in \overline S_0$.  D'apr\`es le lemme 2, il existe une
hypersurface alg\'ebrique $H=\{P=0\}$ passant par $z$ et 
ne rencontrant pas $\overline
S_0$, o\`u $P$ est un polyn\^ome. Pour $\epsilon>0$ suffisamment
petit, on a $bS_0\cap P^{-1}(|x|\leq \epsilon)=\emptyset$. 
Par cons\'equent, $\frac{1}{2\pi
  i}\int_{bS_0}\frac{dP}{P-\alpha}=0$ pour tout
$|\alpha|\leq\epsilon$. Pour tout $i$ suffisamment grand, 
l'int\'egrale  $\frac{1}{2\pi
  i}\int_{bS_i}\frac{dP}{P-\alpha}$ existe car $bS_i$ tend vers $bS_0$
au sens g\'eom\'etrique. Cette int\'egrale est \'egale au nombre de points
d'intersection de $H_\alpha:=P^{-1}(\alpha)$ avec $S_i$. 
Ce nombre entier tend vers 0 car $bS_i$ tend vers $bS_0$ au sens des
courants. Par
cons\'equent, il est \'egal \`a 0 pour $i$ assez grand. Ceci montre
que $z$ n'appartient pas \`a la limite de $\overline S_i$ car $z$ est
un point d'int\'erieur de $P^{-1}(|x|\leq\epsilon)$. 
\end{preuve}
\begin{lemme} Soient $K$ un compact m\'eromorphiquement convexe dans
 une vari\'et\'e de Stein $Z$ et $H$ un compact de mesure de
 Haussdorff $2$-dimensionnelle ${\cal H}^2$ nulle. Alors $K\cup H$
 est m\'eromorphiquement convexe dans $Z$.
\end{lemme}
\begin{preuve}
Comme $Z$ est de Stein, 
il suffit de consid\'erer $Z=\mathbb{C}^n$. Soit $z\notin K\cup
H$. Comme $K$ est m\'eromorphiquement convexe, il existe un polyn\^ome
$P$ v\'erifiant $z\in P^{-1}(0)\subset\mathbb{C}^n\setminus K$. Comme
${\cal H}^2(H)=0$, il existe un polyn\^ome $Q$ v\'erifiant $Q(z)=0$ et
$Q^{-1}(0)\cap P^{-1}(0)\cap H=\emptyset$. Posons $T=Q^{-1}(0)\cap
P^{-1}(0)$, la famille d'hypersurfaces $\{P+\alpha
Q\}_{\alpha\in\mathbb{C}}$ forme un feuilletage holomorphe de
$\mathbb{C}^n\setminus T$. Comme ${\cal
  H}^2(H)=0$, il existe un $\alpha$ assez petit tel que $\{P+\alpha
Q=0\}\cap H=\emptyset$. Pour $\alpha$ assez petit $\{P+\alpha
Q=0\}\cap K=\emptyset$. Alors $\{P+\alpha Q=0\}$ est une hypersurface
passant par $z$ et elle ne coupe pas $K\cup H$. 
D'o\`u $K\cup H$ est m\'eromorphiquement
convexe.
\end{preuve}
\begin{lemme}
Soient $U,V\subset Z$ deux ouverts connexes de Stein, $K\subset U$ un
compact et $H\subset V$ un compact connexe de
mesure ${\cal H}^2$ nulle. Alors $K\cup H$ admet un voisinage de
Stein.
\end{lemme}
\begin{preuve}
Comme $U$ et $V$ sont des ouverts connexes 
de Stein, ils admettent des suites 
d'exhaustion de compacts connexes
holomorphiquement convexes. Par cons\'equent, il existe un compact
connexe, holomorphiquement convexe $K'\subset U$, contenant $K$. Soit $U_1$ un
ouvert connexe \`a bord $C^\omega$ strictement pseudoconvexe de $U$ et
contenant $K'$. 
Si $K'\cap H=\emptyset$ il suffit de choisir un
  voisinage de Stein $W_1$ de $K'$ dans $U_1\setminus H$ et 
gr\^ace au lemme 3, un voisinage de Stein
  $W_2$ de $H$ dans $V\setminus K'$. L'ouvert 
 $W=W_1\cup W_2$ v\'erifie notre lemme. 
Sinon, $F_1:=H\cap bU_1\not =\emptyset$.
  Soit $U_2$ un autre ouvert \`a bord $C^\omega$ dans $U$, contenant
  $\overline U_1$. On a \'egalement  $F_2:=H\cap bU_2\not
  =\emptyset$. Dans $V$, $F_1$ est m\'eromorphiquement convexe, il
  admet un syst\`eme fondamental de voisinages m\'eromorphiquement
  convexes. Soit $G$ un
  voisinage m\'eromorphiquement convexe 
de $F_1$ dans $V$ v\'erifiant $G\cap K'=\emptyset$ et
  $\overline G\subset U_2$. D'apr\`es le lemme pr\'ec\'edent, $(K'\cup
  H)\cap U_2$ est m\'eromorphiquement convexe dans $U$. Il admet donc
  un syst\`eme fondamental de voisinages de Stein. On choisit
  $W_1:=\{\rho_1<0\}$ un voisinage de Stein v\'erifiant $W_1\cap
  bU_1\subset G$ o\`u $\rho$ est une fonction ${\cal C}^\omega$
  strictement p.s.h. d\'efinie au voisinage de $W_1$. 
D'apr\`es le lemme 3, $\overline G\cup (H\setminus U_1)$ est
m\'eromorphiquement convexe dans $V$. 
On choisit $W_2:=\{\rho_2<0\}$ un voisinage de
  Stein de $\overline G\cup (H\setminus U_1)$ v\'erifiant $W_2\cap
  bU_2\subset W_1$, o\`u $\rho_2$ est une fonction ${\cal C}^\omega$ 
strictement p.s.h
  au voisinage de $W_2$. On d\'efinit un voisinage ouvert $W^*$ de
  $K'\cup H$ par 
\begin{enumerate}
\item $W^*\cap U_1:=\{\rho_1<0\}\cap U_1$.
\item $W^*\cap (U_2\setminus U_1):=\{\rho_1<0,\rho_2<0\}\cap (U_2\setminus
  U_1)$.
\item $W^*\setminus U_2:=\{\rho_2<0\}\setminus U_2$.
\end{enumerate}
Soit $W^{**}$ la composante connexe de $W^*$ contenant $K'\cup H$.
Sur son bord, le domaine $W^{**}$ est d\'efini localement par une ou deux
in\'egalit\'es $\varphi<0$ o\`u $\varphi$ est une fonction strictement
p.s.h. Si $W^{**}$ ne contient aucun sous-ensemble analytique compact
de dimension $\geq 1$,
d'apr\`es le th\'eor\`eme de Grauert \cite{Grauert}, $W^{**}$ est de
Stein. Sinon, $W^{**}$ contient un  sous-ensemble
analytique compact maximal
$L$,  dont toute composante irr\'eductible 
est de dimension $\geq 1$. On appelle
$L_1,\ldots, L_m$ les composantes irr\'eductibles de $L$.
Comme $U$ est de Stein et $H$ est un compact de $V$, on a $L_i\cap
V\setminus H\not=\emptyset$ pour tout $i$. Dans $V$, d'apr\`es le
lemme 3, $H$ est m\'eromorphiquement convexe. Par cons\'equent, il
existe une hypersurface $H_i$ de $V$ telle que $H_i\cap H=\emptyset$,
$H_i\cap L_i\not=\emptyset$. Alors $V':=V\setminus\bigcup_{i=1}^m H_i$
est un ouvert de Stein. Soit $U''$ (resp $V''$) un ouvert de $U$
(resp. de $V'$) contenant $K'$ (resp. $H$), \`a bord $C^\omega$,
strictement pseudoconvexe.
Posons $U^*:=U''\cap W^{**}$ et
${V'}^*:=V''\cap W^{**}$. Ces ouverts sont inclus dans des ouverts de Stein
et \`a bord strictement pseudoconvexe. Ils sont donc de Stein. On fait
la m\^eme construction comme celle pr\'ec\'edente, en rempla\c cant
$U$ par $U^*$ et $V$ par $V^*$, pour obtenir un
voisinage connexe $W$ \`a bord strictement pseudoconvexe de $K\cup H$. Cet
ouvert $W$ est inclus dans $W^{**}$ et il ne contient aucun
sous-ensemble analytique 
compact de dimension $\geq 1$ de $W$. 
Il est, d'apr\`es Grauert, ouvert de Stein.
\end{preuve}
\begin{lemme} Toute r\'eunion finie $\Gamma$ d'arcs r\'eels compacts de
  longueur finie d'une vari\'et\'e complexe $Z$ admet un voisinage
  connexe de Stein. Par cons\'equent, $\Gamma$ admet un syst\`eme
  fondamentale de voisinages de Stein et pour tout ensemble fini
  $A\subset Z$, il existe un disque holomorphe, \`a bord lisse dans $Z$
  contenant $A$.
\end{lemme}
\begin{preuve} On choisit un arc r\'eel, ferm\'e de longueur finie $L$
  dans $Z$ tel que $\Gamma':=\Gamma\cup L$ soit connexe. 
On recouvre $\Gamma'$ par une famille finie d'ouverts isomorphe \`a la
boule unit\'e de $\mathbb{C}^n$, o\`u $n=\dim Z$. On note
$U_1,U_2,\ldots, U_m$ ces ouverts. Il existe des compacts $\Gamma_k\subset
\Gamma'\cap U_k$ tels que $\Gamma'=\bigcup\Gamma_k$. On peut supposer
de plus que $\Gamma_k$ sont tous des r\'eunions finies d'arcs
r\'eels ferm\'es. On choisit $L_k\subset U_k$ des arcs r\'eels
ferm\'es tels que $\Gamma_k':=\Gamma_k\cup L_k$ soient connexes. 
On peut num\'eroter
les $U_i$ de sorte que $\bigcup_{i=1}^k \Gamma_i'$ soit connexe pour
tout $k$. En appliquant le
lemme pr\'ec\'edent \`a $U:=U_1$, $V:=U_2$, $K:=\Gamma_1'$ et
$H:=\Gamma_2'$, on peut construire un voisinage de Stein $W_1$ de
$\Gamma_1'\cup\Gamma_2'$. Comme $\Gamma_1'\cup\Gamma_2'$ est connexe, on
peut supposer que $W_1$ est connexe. D'apr\`es le lemme
pr\'ec\'edent appliqu\'e \`a $U:=W_{k-1}$, $V:=U_{k+1}$,
$K:=\bigcup_{i=1}^k\Gamma_i'$ et $H:=\Gamma_{k+1}$, il
existe un voisinage connexe de Stein $W_k$ de
$\bigcup_{i=1}^{k+1}\Gamma_i'$ pour tout $k=2,3,\ldots,m-1$. D'apr\`es
le lemme 3, dans
$W=W_{m-1}$, $\Gamma$ est m\'eromorphiquement convexe. Il admet donc
un syst\`eme fondamental de voisinages de Stein. \\
On choisit  un arc r\'eel ferm\'e de longueur finie $\Gamma^*$ passant
par tous les points de $A$ et un voisinage $U$  connexe et de Stein de
$L$. Soit $S$ un
  sous-ensemble analytique lisse, irr\'eductible, 
de dimension 1
  de $U$ et contenant $A$. Tout ouvert simplement connexe
  \`a bord lisse de  $S$, contenant $A$ v\'erifie notre lemme.
\end{preuve} 
\begin{lemme} Sous l'hypoth\`ese du th\'eor\`eme 2, 
il existe un ouvert connexe $Y$ de $Z$, contenant $S\cup\Gamma$
  tel que dans $Y$, toute surface de Riemann $S'$ 
born\'ee dans $Y$ et v\'erifiant
\begin{enumerate}
\item $S'$ est un sous-ensemble analytique  de $Y\setminus\Gamma$.
\item $\overline S'$ n'est pas un sous-ensemble analytique de $Y$.
\end{enumerate}
soit incluse dans $ S$.
\end{lemme}
\begin{preuve}
Soit $S'$ un sous-ensemble analytique de $Z\setminus\Gamma$, born\'ee
dans $Z$. Dans un voisinage connexe de Stein $V$ de $\Gamma$,
la surface $S'\cap V$ d\'efinit un courant d'int\'egration \`a bord
rectifiable \cite{Lawrence1, Dinh2}. Alors si $d[S'\cap V]=0$,
d'apr\`es le th\'eor\`eme de structure de
King-Harvey-Shiffman-Alexander \cite{King, HarveyShiffmann,Shiffmann,
  Alexander2} $\overline S$ est un sous-ensemble analytique de
$V$. Sinon, $d[S'\cap V]$ est un courant rectifiable, ferm\'e \`a
support dans $\Gamma$. Par cons\'equent,
\begin{enumerate}
\item Soit $\overline S'$ est un sous-ensemble analytique de $Z$.
\item Soit $bS'$ se constitue par des composantes connexes de $\Gamma$. 
\end{enumerate}
Il existe un
nombre fini de surfaces de Riemann de deuxi\`eme type. 
On note $S_1,S_2,\ldots, S_m$ ces surfaces de Riemann dont
$S_k\not\subset S$ pour tout $k<m'$ et $S_k\subset S$ pour tout $k\geq
m'$. Tout ouvert connexe $Y$ 
 contenant $S\cup\Gamma$ tel que $Y\not\supset S_k$ pour tout
 $k<m'$, v\'erifie notre lemme. 
\end{preuve}
On choisit $L$ un arc r\'eel, compact de longueur finie de $Y$ tel que
$\Gamma\cup L$ soit connexe, $L\cap S=\emptyset$ et tel que le lemme
pr\'ec\'edent reste valable si l'on remplace $\Gamma$ par
$\Gamma':=\Gamma\cup L$.\\
D'apr\`es le lemme 5, il existe un voisinage connexe de Stein
$U\subset Y$
de $\Gamma'$. On choisit un ouvert connexe $U_1\subset U$ 
strictement pseudoconvexe, \`a bord
lisse et contenant $\Gamma'$ tel que $bU_1$ coupe $S$
transversalement en un nombre fini de courbes r\'eelles compactes, lisses
$C_1,\ldots, C_m$ de $S$. D'apr\`es le lemme 2, dans $U$,
$\widehat{C}\setminus C$ est un sous-ensemble analytique
de dimension 1 de $U\setminus C$, born\'e dans $U$, o\`u
$C:=\bigcup_{k=1}^m C_k$. D'apr\`es le lemme pr\'ec\'edent,
$\widehat{C}\setminus C\subset \overline S$.
\begin{trivlist}
\item[] \textit{Premier cas.} Supposons que $\widehat{C}\setminus
  C=\emptyset$. Alors $C$ est holomorphiquement convexe dans
  $U$. Soit $U_2\subset U$ (resp. $U_0$) un ouvert \`a bord lisse
  contenant $\overline U_1$
  (resp. relativement compact dans $U_1$),
  qui est une petite d\'eformation lisse de $U_1$. Alors $bU_2\cap S$
  (resp. $bU_0\cap S$)
  est encore holomorphiquement convexe dans $U$. On choisit $G\subset\subset
  U\setminus U_1$ un voisinage de $bU_2\cap S$ tel que
  $\overline G$ soit holomorphiquement convexe dans $U$. D'apr\`es le
  lemme 1, $\widehat{\overline G\cup\Gamma'}\setminus \overline
  G\cup\Gamma'$ est un sous-ensemble analytique de dimension 1 de
  $U\setminus (\overline G\cup\Gamma')$. Par cons\'equent, 
d'apr\`es le lemme 6, si $G$
  est un voisinage suffisament petit de $bU_2\cap S$, on a
  $\widehat{\overline G\cup\Gamma'} 
\setminus (\overline G\cup\Gamma')\subset S\cup\Gamma'\cup (Y\setminus
  \overline U_1)$. D'apr\`es le
  th\'eor\`eme de Siu, il existe un voisinage de Stein $W_1$ 
  de $S\setminus \overline U_0$ tel que $W_1\cap bU_2\subset
  G$. Soit $W_2=\{\rho_2<0\}\subset W_1$ un ouvert \`a bord lisse,
  strictement pseudoconvexe et contenant $S\setminus U_1$, o\`u
  $\rho_2$ est une fonction strictement p.s.h et d\'efinie au
  voisinage de $W_2$. Soit $W_3:=\{\rho_3<0\}$ un voisinage de Stein
  de $\widehat{\overline G\cup\Gamma'}$ v\'erifiant $W_3\cap
  bU_1\subset W_2$,  o\`u
  $\rho_3$ est une fonction strictement p.s.h et d\'efinie au
  voisinage de $W_3$. On d\'efinit le voisinage $W^*$ de $S\cup\Gamma'$
  par
\begin{enumerate}
\item $W^*\setminus U_2=\{\rho_2<0\}\setminus U_2$.
\item $W^*\cap (U_2\setminus U_1)=\{\rho_2<0, \rho_3<0\}\cap
  (U_2\setminus U_1)$.
\item $W^*\cap U_1=\{\rho_3<0\}\cap U_1$.
\end{enumerate}
Soit $W^{**}$ la composante connexe de $W^*$ contenant $S\cup\Gamma'$.
Sur son bord, cet ouvert est localement d\'efini  par une ou deux
in\'egalit\'es $\varphi<0$ o\`u $\varphi$ est une fonction strictement
p.s.h. Si $W^{**}$ ne contient aucun sous-ensemble analytique compact
de dimension $\geq 1$, 
d'apr\`es Grauert, il est de Stein \cite{Grauert}. 
Sinon, $W^{**}$ contient un sous-ensemble analytique compact maximal
dont chaque composante irr\'eductible  est de dimension $\geq 1$. On
appelle $H_1,\ldots, H_m$ ces composantes irr\'eductibles.
 Par
hypoth\`ese, aucune $H_i$
n'est incluse dans $S\cup\Gamma'$. On remplace $Y$ par un ouvert
connexe de
$W^{**}$ qui contient $S\cup \Gamma'$ et ne contient aucune $H_i$. La
m\^eme construction que celle pr\'ec\'edente nous donne un voisinage
connexe $W$ de $S\cup\Gamma'$ dont le bord est strictement
pseudoconvexe. L'ouvert $W$ ne contient aucun sous-ensemble analytique
compact de dimension $\geq 1$.
Il est donc ouvert de Stein.  D'apr\`es le lemme 2, $S\cup\Gamma$
admet un syst\`eme fondamental de voisinages de Stein.
\item[] \textit{Deuxi\`eme cas.} Supposons que $(\widehat C\setminus
  C)\cap \Gamma'=\emptyset$. Posons $T:=\widehat C\setminus
  C$. C'est un sous-ensemble analytique de dimension 1 de $U_1$. Soit $H$
  une hypersurface de $U_1$, contenant $T$
  v\'erifiant $H\cap\Gamma'=\emptyset$ et $\#H\cap \overline{S\setminus
  T}<+\infty$. Alors 
$U_1\setminus H$ est de Stein. Il poss\`ede une suite d'exhaustion
  d'ouverts connexes $\C^\omega$ 
strictement pseudoconvexes $V_k$. Pour $k$ assez grand,
  $bV_k\cap S$ est constitu\'e par 
  un nombre fini de petites courbes ferm\'ees qui bordent un petit voisinage
  dans $\overline{S\setminus
  T}$ de l'ensemble fini $H\cap \overline{S\setminus
  T}$. Par cons\'equent, comme $S$ et $Y$ v\'erifie le lemme 6 et
  comme $\overline S$ ne contient aucune surface de Riemann compacte, 
$\tilde C:=bV_k\cap S$ est holomorphiquement
  convexe dans $V_{k+1}$. En rempla\c cant $U$ par $V_{k+1}$ et $U_1$
  par $V_k$, on se ram\`ene au premier cas.
\item[] \textit{Cas g\'en\'eral.} D'apr\`es le th\'eor\`eme
  d'unicit\'e,  $M:=(\widehat C\setminus C)\cap\Gamma'$ est 
un ferm\'e de mesure
  ${\cal H}^1$ nulle. Par cons\'equent, d'apr\`es le lemme 1, 
il est holomorphiquement
  convexe dans $U$. 
\begin{lemme}
Il existe les voisinages 
  holomorphiquement convexes, \`a
  bord lisse $V_0$, $V_1$, $V_2$ de $M$ v\'erifiant les conditions suivantes
\begin{enumerate}
\item $V_{i+1}$ est un ouvert relativement compact de $V_i$.
\item $V_i$ poss\`ede un nombre fini de composantes connexes.
\item Le bord de $V_i$ intersecte $S\cup\Gamma'$ transversalement et
  $bV_i\cap \Gamma'$ est un ensemble fini.
\item $(S\cup\Gamma')\cap bV_1$ est holomorphiquement convexe dans
  un voisinage de Stein de $(S\cup\Gamma')\setminus V_2$.
\end{enumerate}
\end{lemme}
\begin{preuve}
Supposons que le lemme est faux. On peut
facilement construire  une suite 
d'ouverts holomorphiquement convexes $\{V_i\}_{i=0}^\infty$ 
d\'ecroissante vers $M$
et
v\'erifiant les conditions 1-3 pour $i=0,1,2,\ldots$. 
Posons $\Gamma^*:=(S\cap bV_i)\cup (\Gamma'\setminus V_i)$ et
$S^*:=S\setminus V_i$. D'apr\`es le lemme 1, $(S^*\cup\Gamma')\cap
\overline U_1$ est m\'eromorphiquement convexe dans $U$. En plus, dans
$U$, $\widehat{\tilde\Gamma}\setminus{\tilde\Gamma}$ se constitue par un
nombre fini de surfaces de Riemann irr\'eductibles, o\`u $\tilde
\Gamma:= (S^*\cap b\overline U_1)\cup(\Gamma'\cap\overline U_1)$. Par
cons\'equent, il existe une hypersurface $H$ de $U$ telle que
$(S^*\cup\Gamma')\cap\overline U_1$ est holomorphiquement convexe dans
l'ouvert de Stein $U\setminus H$. Alors $H\cap S\cap V_i$ est un
ensemble fini. Comme $U_1\setminus H$ est de Stein, il admet une
suite d'exhaustion d'ouverts connexes holomorphiquement convexes
$\Omega_k$. Pour $k$ suffisament grand, $C^*:=b\Omega_k\cap S$ se
constitue par une petite d\'eformation de $C$ et par des courbes de
Jordan qui forment le bord d'un petit voisinage de $H\cap S\cap V_i$
dans $S$. D'apr\`es le deuxi\`eme cas (appliqu\'e \`a
$U^*:=\Omega_{k+1}$, $U_1^*:=\Omega_k$, $S^*$, $\Gamma^*$ et $C^*$) 
 $S^*\cup\Gamma^*$ admet un syst\`eme
fondamental de voisinages de Stein. Par cons\'equent, comme la
condition 4 n'est pas v\'erifi\'ee quand on remplace les ouverts $V_1$
et $V_2$ dans le lemme par les ouverts $V_{i-1}$ et 
$V_{i}$,
$(S\cup\Gamma')\cap bV_{i-1}$ doit contenir le bord d'une surface de Riemann
$S_i\subset \overline S$. 
Comme $V_i$ est holomorphiquement convexe dans $U$, la surface
$S_i$ ne peut pas \^etre enti\`erement dans $U$. 
Rappellons que quand $i$ tend vers infini, $V_i$
tend vers $M$ qui est de longueur $0$. Finalement, 
on peut conclure ici que $\overline S$
contient une surface de Riemann compacte, qui est la limite de $S_i$
quand $i$ tend vers l'infini. C'est une contradiction. 
\end{preuve}
On choisit d'apr\`es les lemmes 7 et 2,  un voisinage de Stein
$W_0\subset Y\setminus\overline V_2$ de $(S\cup \Gamma')\setminus V_1$
dans $Y$ v\'erifiant 7.4. 
Soit $G\subset W_0$ un voisinage holomorphiquement convexe
de   $(S\cup \Gamma')\cap bV_1$ dans $W_0$. Posons
$\Gamma^*:=(S\cap bV_1)\cup (\Gamma'\setminus V_1)$. D'apr\`es le
lemme 1, dans $W_0$, $\widehat{\overline G\cup\Gamma^*}\setminus
(\overline G\cup \Gamma^*)$ est un sous-ensemble analytique de
dimension 1 de $W_0\setminus (\overline G\cup \Gamma^*)$. Rappellons
que $(S\cup\Gamma')\setminus V_1$ est holomorphiquement convexe dans
$W_0$. Par cons\'equent, d'apr\`es les lemmes 6 et 1, 
si $G$ est un voisinage assez petit, 
$\widehat{\overline G\cup\Gamma^*}\setminus
(\overline G\cup \Gamma^*)$ est inclus dans $\overline S\cup V_0$.
D'apr\`es les lemmes 6 et 1, $(S\cup\Gamma')\cap
\overline V_0$ est holomorphiquement convexe dans $U$. Il poss\`ede donc un
voisinage de Stein $W_1:=\{\rho_1<0\}$ v\'erifiant
$\overline W_1\cap b V_1\subset G$ o\`u $\rho_1$ est une fonction
$C^\omega$, strictement p.s.h au voisinage de $\overline W_1$.
Il
  existe un voisinage de Stein  $W_2:=\{\rho_2<0\}$ 
de  $\widehat{\overline G\cup\Gamma^*}$
  v\'erifiant
$\overline W_2\cap b V_0\subset W_1$ o\`u $\rho_2$ est une fonction
$C^\omega$, strictement p.s.h au voisinage de $\overline W_2$. On d\'efinit
le voisinage $W^*$ de $S\cup\Gamma'$ par 
\begin{enumerate}
\item $W^*\cap V_1:=\{\rho_1<0\}\cap V_1$.
\item $W^*\cap (V_0\setminus V_1):=\{\rho_1<0,\rho_2<0\}\cap
  (V_0\setminus V_1)$.
\item $W^*\setminus V_0:=\{\rho_2<0\}\setminus V_0$
\end{enumerate}
Soit $W^{**}$ la composante connexe de $W^*$ contenant $S\cup\Gamma$. 
Cet ouvert est \`a bord strictement
pseudoconvexe, il peut poss\'eder un sous-ensemble analytique compact
maximal dont toute composante irr\'eductible  est de dimension
$\geq 1$. On appelle $H_1,\ldots,H_m$ ces
composantes irr\'eductibles.
En rempla\c cant $Y$ par un ouvert de $W^{**}$ qui ne contient aucune
$H_i$, on peut construire par la m\'ethode
pr\'ec\'edente, un ouvert connexe \`a bord strictement pseudoconvexe
$W\subset W^{**}$ contenant $S\cup \Gamma'$. Cet ouvert ne contient
aucun sous-ensemble analytique compact de dimension $\geq 1$. 
Il est, d'apr\`es Grauert, ouvert de
Stein. D'apr\`es le lemme 2,
$S\cup\Gamma$ poss\`ede dans $W$ un syst\`eme fondamental 
de voisinages de Stein.
\end{trivlist}
\section{Preuve du th\'eor\`eme 1}
D'apr\`es le th\'eor\`eme de Dolbeault-Henkin g\'en\'eralis\'e
\cite{Dinh4}, il suffit de prouver que le plan tangent de 
$\Gamma$ est maximalement
complexe en tout point de $\Gamma\cap
\mathbb{P}^{n-p+1}_\nu$ pour un ensemble
dense de $\nu\in V$. 
Par m\'ethode de tranchage, on ram\`ene le probl\`eme au
cas $p=2$. Par m\'ethode de projection, le probl\`eme se ram\`ene
au cas $n=3$ et $p=2$.
\begin{proposition} Pour tout ouvert $V'\subset V$, il existe
  $\nu_0\in V'$, un sous-ensemble $V''\subset V'$ et un voisinage de
  Stein $W$ de $\Gamma\cap\mathbb{P}^2_{\nu_0}$ v\'erifiant
\begin{enumerate}
\item Les vecteurs tangents g\'eom\'etriques de $V''$ en $\nu_0$
  engendent le plan tangent de $V$ en $\nu_0$.
\item Pour tout $\nu\in V''$, l'intersection
  $\Gamma\cap\mathbb{P}^2_\nu$ est incluse dans $W$.
\item Pour tout $\nu\in V''$, l'intersection
  $\Gamma\cap\mathbb{P}^2_\nu$ est le bord d'une $1$-cha\^{\i}ne
  holomorphe $S_\nu$ de $\mathbb{P}^2_\nu\setminus\Gamma$ qui
  est relativement compacte dans $W$.
\end{enumerate}
\end{proposition}
\begin{preuve} On peut supposer que $V'$ est suffisamment petit tel que
  pour tous $\nu,\nu'\in V'$
  la combinaison $\gamma_\nu:=\Gamma\cap\mathbb{P}^2_\nu$ soit une
  d\'eformation continue de $\Gamma\cap\mathbb{P}^2_{\nu'}$. On
  \'ecrit 
$$\gamma_\nu=
\Gamma\cap\mathbb{P}^2_\nu=\sum_{i\in I} n_i\gamma_{i,\nu}$$
o\`u les $n_i$ sont des entiers non nuls 
et $\gamma_{i,\nu}$ sont des courbes
r\'eelles ferm\'ees, irr\'eductibles, $C^2$ d\'ependantes
contin\^ument de $\nu\in V'$. On appelle $(w_0:w_1:w_2:w_3)$ les
coordonn\'ees homog\`enes de $\mathbb{CP}^3$ et $z_i:=w_i/w_0$ ses
coordonn\'ees affines. On peut supposer que les
hyperplans
$H:=\{w_2=0\}$ et $Q=\{w_0=0\}$ ne rencontrent pas $\gamma_\nu$ 
pour tout $\nu\in V'$. De plus, on peut
supposer que $(w_0:w_1:w_2)$ est un syst\`eme des coordonn\'ees de
$\mathbb{P}^2_\nu$  pour
tout $\nu\in V'$. On choisit
pour tout $\nu\in V'$, une $1$-cha\^{\i}ne holomorphe $S_\nu$ de
$\mathbb{P}^2_\nu\setminus\Gamma$ \`a bord $\gamma_\nu$ telle
que $\overline S_\nu$ ne contient aucune surface de Riemann
compacte. Rappelons que deux cha\^{\i}nes de m\^eme bord au sens des courants
diff\`erent d'une cha\^{\i}ne alg\'ebrique. On \'ecrit 
$$S_\nu=\sum_{j\in J_\nu}m_{j,\nu}S_{j,\nu}$$
o\`u  les $m_{j,\nu}$ sont des entiers et $S_{i,\nu}$ sont des
sous-ensembles analytiques irr\'eductibles de
$\mathbb{P}^2_\nu\setminus\Gamma$. Alors le bord de $S_{j,\nu}$ est
une combinaison lin\'eaire de $\gamma_{j,\nu}$ \`a coefficients $0$,
$1$ ou $-1$. Soit $F_\nu$ l'ensemble de toutes ces combinaisons. Comme
$V'$ est petit et $\gamma_\nu$ d\'epend contin\^ument de $\nu$, on
peut identifier tous les
$F_\nu$ et noter $F:=F_\nu$ pour $\nu\in V'$. On peut \'egalement
consid\'erer $J_\nu$ comme un sous-ensemble de $F$. On appelle
$k_{j,\nu}$ le nombre de point d'intersection de $S_{j,\nu}$ avec $H$
en comptant les multiplicit\'es. Soient $J$ un sous-ensemble de $F$,
$M=\{m_j\}_{j\in J}$ et $K=\{k_j\}_{j\in J}$ des familles de nombres
entiers. On pose
$$V_{J,M,K}:=\{\nu\in V'' \mbox{ telle que } J_\nu=J, m_{j,\nu}=m_j,
k_{j,\nu}=k_j \mbox{ pour tout } j\in J\}.$$
Alors $V'=\bigcup V_{J,M,K}$. Par cons\'equent, 
il existe un $V_{J,M,K}$ et un $\nu_0\in V_{J,M,K}$ tels que les
vecteurs tangents g\'eom\'etriques de $V_{J,M,K}$ en $\nu_0$
engendent le plan tangent de $V$ en $\nu_0$. D'apr\`es le
th\'eor\`eme 2, il suffit de prouver que $S_{j,\nu}$ d\'epend de
mani\`ere continue de $\nu\in V_{J,M,K}$. Fixons un $j\in J$. Pour
simplifier des notations, on pose $k:=k_j$.
Soient $X'$ un voisinage 2-concave suffisamment petit de 
$\mathbb{P}^2_{\nu_0}$ et $\Gamma'$ une
combinaison lin\'eaire \`a coefficients $\pm 1$ de composantes
irr\'eductibles de $\Gamma\cap X'$ tels que $\gamma'_\nu:=\Gamma'\cap
\mathbb{P}^2_\nu= bS_{j,\nu}$ pour $\nu\in V_{J,M,K}$ proche de
$\nu_0$. On pose 
$$G_{i,\nu}(\xi,\eta):=
\frac{1}{2\pi i}\int_{\gamma'_\nu}z_1^i\frac{d(z_2-\xi-\eta
       z_1)}{z_2-\xi-\eta z_1}.$$
Ces fonctions sont les fonctions sym\'etriques qui permettent de
       d\'eterminer les solutions du probl\`eme du bord pour
       $\gamma'_\nu$ \cite{DolbeaultHenkin}.
Un th\'eor\`eme de Dolbeault-Henkin implique que les deux
conditions suivantes sont \'equivalentes ({\it voir} \cite{Dinh4})
\begin{enumerate}
\item $\gamma'_\nu$ est le bord d'une $1$-cha\^{\i}ne holomorphe $S'_\nu$
  telle que l'intersection de $H$ avec $S'_\nu$ soit une combinaison
  de points \`a coefficients entiers positifs et de masse totale $k$.
\item Il existe des fonctions $C_{i,\nu}(\eta)$ d\'efinies $C^\infty$ dans
 un ouvert $\Omega$ de $\mathbb{C}$ pour $i=0,1,\ldots, k+1$ et un voisinage
 $U$ de $0$
 dans $\mathbb{C}$ tels que pour tout $\xi\in U$ 
le syst\`eme d'\'equations suivant admette une solution
$$\left\{ \begin{array}{c}\displaystyle
 x_1^i+\cdots +x_{k}^i
     = G_{i,\nu}(\xi,\eta)+\frac{C_{0,\nu}^{(i-1)}(\eta)}{(i-1)!}\xi^i
 +\displaystyle\sum_{s=1}^i
\frac{iC_{s,\nu}^{(i-s)}(\eta)}{(i-s)!}\xi^{i-s}\\
\\
\mbox{ pour } i=1,2,\ldots, k+1
\end{array}
\right. $$
\end{enumerate}
Si, dans la condition 1, la solution $S'_\nu$ existe, 
elle sera unique. Par cons\'equent,
pour $\nu\in V_{J,M,K}$, cette solution est \'egale \`a
$S_{j,\nu}$. Au voisinage de $H\cap\mathbb{P}^2_\nu$, $S'_\nu$ sera
d\'ecrite par $\{(x_s,\xi+\eta x_s)\}_{s=1}^k$ pour $(\xi,\eta)$
variant dans $U\times\Omega$ \cite{Dinh4}. Soient $S_i:=x_1^i+\cdots
+x_{k}^i$. Alors $S_{k+1}$ sera d\'etermin\'ee par un polyn\^ome en
$S_1$, ..., $S_k$. On \'ecrit
$$S_{k+1}=P(S_1,\ldots,S_k).$$
On remplace les $S_i$ dans l'\'equation pr\'ec\'edente par les membres
\`a droite du syst\`eme ci-dessus. Ensuite, on d\'eveloppe deux
membres de l'\'equation obtenue en s\'erie de Taylor en
$\xi$. L'identification des coefficients de $\xi^i$ nous donne un
syst\`eme de la forme
$$\left\{\begin{array}{c}
\displaystyle C_{k,\nu}^{(k+1-i)}(\eta)=F_i(\eta,C_{0,\nu} ,
C_{0,\nu}^{(1)},\ldots, C_{0,\nu}^{(k-1)}, C_{1,\nu}, C_{1,\nu}^{(1)},
\ldots, C_{k,\nu},\nu)\\
\mbox{ pour } i=k+1,k,\ldots,1
\\
C_{0,\nu}^{(k)}(\eta)=F_0(\eta,C_{0,\nu} ,
C_{0,\nu}^{(1)},\ldots, C_{0,\nu}^{(k-1)}, C_{1,\nu}, C_{1,\nu}^{(1)},
\ldots, C_{k,\nu},\nu)
\\
H_s(\eta,C_{0,\nu},\ldots,C_{0,\nu}^{(k-1)},\ldots,C_{k,\nu},\nu)=0\\

\mbox{ pour } s=1,2,\ldots
\end{array}
\right. $$
où $F_i$ et $H_s$ 
sont des polyn\^omes en $C_{s,\nu}^{(i)}$, holomorphes en $\eta$ et
$C^2$ en $\nu$.\\
D'apr\`es le th\'eor\`eme de Cauchy-Kovalevskaya, 
les $k+1$ premi\`eres \'equations donnent une solution unique pour
chaque donn\'ee vectorielle  
$b$ en un point fix\'e $\theta$ de $\Omega$, c.-\`a-d. 
$$(C_{\nu,0},\ldots,C_{0,\nu}^{(k-1)}, C_{1,\nu}, C_{1,\nu}^{(1)},
\ldots, C_{k,\nu})=b \mbox{ en }\theta.$$  
Cette solution
d\'epend holomorphiquement de $(b,\eta)$ et de fa\c con $C^2$ de
$\nu$. Rempla\c cant cette solution dans le reste du syst\`eme, nous
obtenons un syst\`eme de la forme
$$\left\{\begin{array}{c}
\displaystyle
\tilde H_s(b,\eta,\nu)=0\\
\mbox{ pour } s=1,2,\ldots
\end{array}
\right. $$
o\`u $\tilde H_s$ sont des fonctions holomorphes en $b$, $\nu$ et
$C^2$ en $\nu$.\\
L'ensemble de solutions de ce syst\`eme repr\'esente un ferm\'e
$\Sigma$ dans
l'espace de $(b,\eta,\nu)$. Le fait
que  pour
un $\nu$, il existe une $1$-cha\^{\i}ne $S'_\nu$ v\'erifiant la
condition 1 ci-dessus, est \'equivalent au fait qu'il existe un $b$ et
un ouvert non vide $\Theta$ de $\mathbb{C}$ v\'erifiant
$\{b\}\times\Theta\times \{\nu\}\subset \Sigma$. 
Ceci montre que la solution du syst\`eme
initial d\'epend de mani\`ere continue de $\nu\in V_{J,M,K}$. 
Soit
$X''\subset X'$ un petit voisinage suffisamment petit 
de $H\cap \mathbb{P}^2_{\nu_0}$, \`a bord
lisse tel que $bX''$ coupe $S'_{\nu_0}$ transversalement.  
Alors $S'_\nu\setminus  X''$ tend vers $S'_{\nu_0}\setminus
X''$ quand $\nu\in V_{J,M,K}$ tend vers $\nu_0$. Posons
$S''_\nu:=S'_\nu\setminus \overline X''$. D'apr\`es le
corollaire 1, 
$S''_\nu$ tend vers $S''_{\nu_0}$ 
quand $\nu\in V_{J,M,K}$ tend vers $\nu_0$. D'o\`u $S_{j,\nu}$ tend
vers $S_{j,\nu_0}$ quand $\nu\in V_{J,M,K}$ tend vers $\nu_0$.
\end{preuve}
La preuve du th\'eor\`eme 1 sera compl\'et\'ee par la proposition
suivante:
\begin{proposition} Dans les conditions du th\'eor\`eme 1 et de la
  proposition 1, le plan tangent de $\Gamma$ en chaque point de
  $\Gamma\cap\mathbb{P}^2_{\nu_0}$ est maximalement complexe. 
\end{proposition}
\begin{preuve} Soient $\Pi:\
  \mathbb{CP}^3\setminus (0:0:0:1)\longrightarrow
  \mathbb{CP}^2$ la projection canonique $\Pi(w):=(w_0:w_1:w_2)$ et
  $\Pi_1:\ \mathbb{CP}^3\setminus \{w_0=0\}\longrightarrow
  \mathbb{C}$ avec $\Pi_1(z):=z_1$. 
On peut supposer que $\{w_0=0\}\cap
  \gamma_\nu=\emptyset$ pour tout $\nu$ proche de $\nu_0$,
  $(0:0:0:1)\in\mathbb{P}^2_{\nu_0}$, la restriction de $\Pi_1$ sur
  $\gamma_{\nu_0}$ est $C^2$, injective en dehors d'un ensemble fini et que la
  restriction de la projection $\Pi$ sur $\Gamma\cap U$ soit de rang maximal
  en tout point et injective en dehors d'une hypersurface r\'eelle de
  $\Gamma\cap U$, o\`u $U$ est un petit voisinage de
  $\gamma_{\nu_0}$. La vari\'et\'e $\Gamma\cap U$ peut \^etre
  consid\'er\'ee comme le graphe d'une fonction $Z_3$ au dessus de
  $\underline\Gamma:=\Pi(\Gamma)\cap \underline U$ sauf sur les
  singularit\'es 
de $\underline\Gamma$, o\`u $\underline U$ est un petit voisinage de
  $\Pi(\mathbb{P}^2_{\nu_0})\simeq\mathbb{CP}^1$. La
  fonction $Z_3$ est $C^2$ sur son domaine de d\'efinition. Comme $W$
  est un ouvert de Stein, il existe une fonction $g$ non identiquement
  nulle et holomorphe dans $W$ telle que $g^{-1}(0)\supset
  W\cap\{w_0=0\}$. D'apr\`es la proposition 1, la
  condition des moments (ou la formule de Stokes sur $S_\nu$) implique que
  $\int_{\gamma_\nu}P(z_1,z_2)g^kdz_1=0$ pour tout  
polyn\^ome $P$, tout $k\geq\deg P+2$ et
  pout tout $\nu\in V''$. \\
On appelle $L$ l'hyperplan complexe de l'espace
  $\mathbb{C}\otimes_\mathbb{R}\Tan (\G(3,4),\nu_0)$, engendr\'e par les
  vecteurs tangents de $V''$ en $\nu_0$, o\`u la notation
  $\Tan(\G(3,4),\nu_0)$
  signifie le plan tangent r\'eel de la grassmannienne $\G(3,4)$ 
en $\nu_0$. Alors pour tout vecteur
  ${\cal L}\in L$ la d\'eriv\'ee  
${\cal L}(\int_{\gamma_\nu}P(z_1,z_2)g^kdz_1)$ s'annulle en  $\nu_0$ pour
  $k\geq \deg Q+2$.\\
Soit $M:=\{\nu\in\G(3,4):\ \mathbb{P}^2_\nu\ni (0:0:0:1)\}$. Posons
  $\Lambda:= L\cap\mathbb{C}\otimes_\mathbb{R}\Tan(M,\nu_0)$. Alors
  $\Lambda$ est un hyperplan de
  $\mathbb{C}\otimes_\mathbb{R}\Tan(M,\nu_0)$. On peut identifier $M$
  \`a la grassmannienne $\G(2,3)$ qui param\`etre l'ensemble de
  droites projectives de $\mathbb{P}^2$. Chaque $\mathbb{P}^2_\nu$
  avec $\nu\in M$ sera associ\'e \`a $\mathbb{P}^1_\nu:=\Pi(\mathbb{P}^2_\nu)$. Alors
  $\dim_\mathbb{C}\Lambda=3$. Une droite projective g\'en\'erique
  $\mathbb{P}^1_\nu$ de $\mathbb{CP}^2$ (avec $\nu\in M$) est
  d\'efinie par l'\'equation $z_2=\xi+\eta z_1$ avec
  $\xi,\eta\in\mathbb{C}$. Alors il existe
  $(a,b), (c,d)\in\mathbb{C}^2\setminus\{0\}$   tels que 
$\Lambda$ est engendr\'e par les trois vecteurs $\mathbb{C}$-ind\'ependants
  $a\frac{\partial}{\partial\xi}+b\frac{\partial}{\partial\eta}$,
  $\overline a\frac{\partial}{\partial\overline \xi}+\overline
  b\frac{\partial}{\partial\overline \eta}$ et
  $c\frac{\partial}{\partial\xi}+\overline
  c\frac{\partial}{\partial\overline \xi}  
+d\frac{\partial}{\partial\eta}+\overline
  d\frac{\partial}{\partial\overline \eta}$. Ces trois vecteurs
  sont $\mathbb{C}$-ind\'ependants si et seulement si la
  matrice suivante est de rang maximal
$${\cal M}:=\left(\begin{array}{cccc}
a & b & 0 & 0\\
0 & 0 & \overline a & \overline b\\
c & d & \overline c &  \overline d
\end{array}\right)$$
Pour tout point r\'egulier 
$z\in\underline\Gamma$, on note $\overline{\cal L}_2$ le vecteur
antiholomorphe de  $\mathbb{C}\otimes_\mathbb{R}
\Tan_\mathbb{C}(\underline\Gamma,z)$ dont la
projection sur $\mathbb{C}\otimes_\mathbb{R} 
\{\frac{\partial}{\partial z_2}\}$
vaut $\frac{\partial}{\partial\overline z_2}$, o\`u la notation 
$\Tan_\mathbb{C}(\underline\Gamma,z)$ signifie la droite tangente
complexe de $\underline\Gamma$ en $z$. Fixons $\nu^*=(\xi^*,\eta^*)\in
M$ et consid\'erons $\nu=(\xi,\eta^*)\in M$. 
Consid\'erons $z^*\in\Pi(\gamma_{\nu_0})$ et $z$ le point d'intersection
de $\mathbb{P}^1_\nu$ avec
$\Tan_\mathbb{C}(\underline\Gamma,z^*)$. Alors la deuxi\`eme
coordonn\'ee $z_2$ de $z$ peut \^etre consid\'er\'ee comme une fonction
holomorphe en $\xi$ et on note $\frac{\partial \overline z_2}{\partial
  \overline \xi}$ pour la d\'eriv\'ee de $\overline z_2$ 
  en $\overline\xi$. Comme le restriction de $\Pi$ sur $\Gamma\cap U$
  est injective en dehors d'une hypersurface r\'eelle de $\Gamma\cap
  U$, on peut poser $\overline{\cal
    L}^*_2(g)(z_1,z_2,z_3):=\overline{\cal
    L}_2(g(z_1,z_2,Z_3(z_1,z_2)))$. C'est une fonction continue sur la
  partie r\'eguli\`ere de $\underline\Gamma$.
\begin{lemme}Avec les notations ci-dessus, pour tout $k\geq 0$ et pour tout
  polyn\^ome $P$ les
  \'egalit\'es suivantes sont vraies en tout point
  $\nu^*=(\xi^*,\eta^*)\in M$
\begin{enumerate}
\item $\displaystyle\frac{\partial \int_{\gamma_\nu}P(z_1,z_2)
g^kdz_1}{\partial
    \eta} =\frac{\partial \int_{\gamma_\nu}P(z_1,z_2)z_1g^kdz_1}{\partial
    \xi}$.
\item $\displaystyle \frac{\partial \int_{\gamma_\nu}P(z_1,z_2)
g^kdz_1}{\partial
    \overline \xi} =
\int_{\gamma_\nu}kP(z_1,z_2)
\overline {\cal L}^*_2 (g)\frac{\partial\overline z_2}
{\partial\overline \xi}g^{k-1}dz_1$.
\item $\displaystyle \frac{\partial \int_{\gamma_\nu}P(z_1,z_2) 
g^kdz_1}{\partial
    \overline \eta} =
\int_{\gamma_\nu}kP(z_1,z_2)\overline z_1
\overline {\cal L}^*_2 (g)\frac{\partial\overline z_2}
{\partial\overline \xi}g^{k-1}dz_1$.
\end{enumerate}
\end{lemme}
\begin{preuve} 
\item {\it Premier cas.} Supposons que $\underline\Gamma$ est une r\'eunion de
droites projectives. On \'ecrit $\underline\Gamma=\bigcup_{\theta\in
  \Upsilon}\mathbb{P}^1_\theta$, o\`u $\Upsilon$ est une combinaison
lin\'eaire \`a coefficients entiers de courbes r\'elles $C^2$ ferm\'ees
dans $M$. La notation $\theta\in \Upsilon$ signifie un point \`a
multiplicit\'e enti\`ere. Pour tout $\theta=(\alpha,\beta)
\in \Upsilon$ et tout $\nu\in M$,
on appelle $z_1(\theta,\nu)$ et $z_2(\theta,\nu)$ les coordonn\'ees du 
point d'intersection de $\mathbb{P}^1_\nu$ avec
$\mathbb{P}^1_\theta$. Alors
\begin{eqnarray} 
\int_{\gamma_\nu}P(z_1,z_2)g^kdz_1 & = & \int_{\theta\in \Upsilon}
\tilde P(\theta,\nu)\tilde g^k(\theta,\nu)dz_1(\theta,\nu)
\end{eqnarray}
o\`u
$$\tilde P(\theta,\nu):=P(z_1(\theta,\nu),z_2(\theta,\nu))$$
  et $$\tilde
  g(\theta,\nu):=g(z_1(\theta,\nu),z_2(\theta,\nu),
Z_3(z_1(\theta,\nu),z_2(\theta,\nu))).$$
Alors $\tilde P$ est holomorphe par rapport \`a $\nu$.
D'autre part on a des relations
$$z_2(\theta,\nu)=\xi+\eta z_1(\theta,\nu)\mbox{\ \ et\ \ } 
z_2(\theta,\nu)=\alpha+\beta z_1(\theta,\nu).$$
Ces relations impliquent
$$\frac{\partial z_1(\theta,\nu)}{\partial\eta}=z_1\frac{\partial
  z_1(\theta,\nu)}{\partial\xi} \mbox{\ \ et\ \ } 
\frac{\partial z_2(\theta,\nu)}{\partial\eta}=z_1\frac{\partial
  z_2(\theta,\nu)}{\partial\xi}$$ 
$$\frac{\partial \overline z_1(\theta,\nu)}{\partial \overline \eta}=
  \overline  z_1\frac{\partial
   \overline z_1(\theta,\nu)}{\partial\xi} \mbox{\ \ et\ \ } 
\frac{\partial  \overline z_2(\theta,\nu)}{\partial \overline \eta}=
  \overline  z_1\frac{\partial
   \overline z_2(\theta,\nu)}{\partial \overline \xi}.$$
Par cons\'equent, pour toute fonction $f$ de classe $C^1$ \`a deux variables
\begin{eqnarray}
\frac{\partial f(z_1(\theta,\nu),z_2(\theta,\nu))}{\partial\eta} 
& = & z_1\frac{\partial
  f(z_1(\theta,\nu),z_2(\theta,\nu))}{\partial\xi}
\end{eqnarray}
et
\begin{eqnarray}
\frac{\partial f(z_1(\theta,\nu),z_2(\theta,\nu))}{\partial\overline\eta} 
& = & \overline z_1\frac{\partial
  f(z_1(\theta,\nu),z_2(\theta,\nu))}{\partial\overline\xi}
\end{eqnarray}
car $z_1(\theta,\nu)$ et $z_2(\theta,\nu)$ sont holomorphes par rapport \`a
$\nu$.\\
Les \'egalit\'es (1), (2) impliquent la premi\`ere \'egalit\'e du
lemme. La deuxi\`eme \'egalit\'e est \'evidente, la troisi\`eme est un
corollaire de (1) et de (3).
\item[] {\it Cas g\'en\'eral.} Pour le cas g\'en\'eral, on consid\`ere
  $\underline\Gamma'$ la r\'eunion des droites tangentes complexes de
  $\underline\Gamma$ aux
  points de $\underline \Gamma'\cap \mathbb{P}^1_{\nu^*}$. 
Alors au voisinage de $\underline \Gamma'\cap \mathbb{P}^1_{\nu^*}$ la
  vari\'et\'e $\underline\Gamma$ est approxim\'ee \`a l'ordre $1$ par
  $\underline\Gamma'$. Ceci explique que les valeurs des d\'eriv\'ees
  d'ordre 1
  trouv\'ees dans le premier cas restent valable dans le cas g\'en\'eral.
\end{preuve}
D'apr\`es le lemme 8, on obtient des \'egalit\'es suivantes en $\nu_0$
pour tout $k\geq \deg P_i+2$
\begin{eqnarray}
\frac{\partial\int_{\gamma_\nu}(a+bz_1)P_1(z_1,z_2)g^kdz_1}{\partial\xi}
& = & 0\\
\int_{\gamma_\nu}(\overline a+\overline b\overline
z_1)P_2(z_1,z_2)\overline {\cal L}^*_2(g)\frac{\partial\overline
  z_2}{\partial\overline \xi} g^{k-1}dz_1
& = & 0\\
\frac{\partial\int_{\gamma_\nu}(c+dz_1)P_3(z_1,z_2)g^kdz_1}{\partial\xi}
+\hspace{3cm}\nonumber \\
+k\int_{\gamma_\nu}(\overline c+\overline d\overline
z_1)P_3(z_1,z_2)\overline {\cal L}^*_2(g)\frac{\partial\overline
  z_2}{\partial\overline \xi} g^{k-1}dz_1
& = & 0
\end{eqnarray}
Pour $P_1=c+dz_1$, $P_2=1$ et $P_3=a+bz_1$, on obtient les
\'egalit\'es suivantes en $\nu_0$ 
pour $k\geq 3$ (la premi\`ere est obtenue par
(5), la seconde par (4) et (6))
\begin{eqnarray}
\int_{\gamma_\nu}(\overline a+\overline b\overline
z_1)\overline {\cal L}^*_2(g)\frac{\partial\overline
  z_2}{\partial\overline \xi} g^{k-1}dz_1
& = & 0\\
\int_{\gamma_\nu}(\overline c+\overline d\overline
z_1)(a+bz_1)\overline {\cal L}^*_2(g)\frac{\partial\overline
  z_2}{\partial\overline \xi} g^{k-1}dz_1
& = & 0
\end{eqnarray}
Comme $W$ est de Stein, il existe une fonction 
$g$ holomorphe dans $W$ v\'erifiant
\begin{enumerate}
\item $g^{-1}(0)\supset W\cap\{w_0=0\}$.
\item L'application $g:\ \gamma_{\nu_0}\longrightarrow \mathbb{C}$
  est $C^2$ et injective en dehors d'un ensemble fini.
\item La d\'eriv\'ee $\partial g/\partial z_3$ ne s'annulle pas sur
  $\gamma_{\nu_0}$.
\end{enumerate}
Si $\overline {\cal L}^*_2(g)$ est identiquement nulle sur $\gamma_{\nu_0}$,
le plan tangent de $\Gamma$ en chaque point de $\gamma_{\nu_0}$ est
maximalement complexe.\\
Supposons que $\overline {\cal L}^*_2(g)$ n'est pas identiquement
nulle sur $\gamma_{\nu_0}$. Posons
$$\gamma:=\overline {\{ z\in\gamma_{\nu_0}:\  {\cal
    L}^*_2(g(z))\not=0 \}}.$$
Alors il existe un arc ouvert $l\subset\gamma_{\nu_0}$ telle que
    $g(l)$ appartient au bord de la composante connexe non
    born\'ee de $\mathbb{C}\setminus g(\gamma)$. Les \'egalit\'es (7)
    et (8) impliquent que les mesures \`a support dans $g(\gamma)$
$$\mu:
=g_*\left\{(\overline a+\overline b\overline
z_1)\overline {\cal L}^*_2(g)\frac{\partial\overline
  z_2}{\partial\overline \xi} g^2dz_1\right\}$$
et 
$$\mu':
=g_*\left\{(\overline c+\overline d\overline
z_1)(a+bz_1)\overline {\cal L}^*_2(g)\frac{\partial\overline
  z_2}{\partial\overline \xi} g^{2}dz_1\right\}$$
sont orthogonales aux polyn\^omes de $\mathbb{C}$. Ces mesures
s'\'etendent dans $\mathbb{C}$ au sens faible des courants en
des $(1,0)$-formes holomorphes dans $\mathbb{C}\setminus
g(\gamma)$ qui sont nulles sur la composante connexe
non born\'ee $R$ de $\mathbb{C}\setminus g(\gamma)$. 
 Sur un petit voisinage d'un 
ouvert dense de $bR$, 
ces extensions sont des extensions continues 
\cite{Dinh3}. On remarque que $S_{\nu_0}$ est lisse en tout point du
bord \`a l'exception d'un compact de longueur $0$. 
Alors il existe un arc ouvert $l'\subset l$ et des ouverts
simplement connexes \`a bord 
lisse $\Theta$, $\Omega$ de $\mathbb{C}$ v\'erifiant
\begin{enumerate}
\item La restriction de l'application $g$ sur $\overline {S_{\nu_0}\cap 
g^{-1}(\Theta)}$ est une application bijective, lisse et
 \`a l'image dans $\overline\Theta$.
\item La restriction de l'application $\Pi_1$ sur $\overline {S_{\nu_0}\cap 
g^{-1}(\Theta)}$ est une application bijective, lisse 
et \`a l'image dans $\overline\Omega$.
\item Les restrictions des mesures $\mu$ et $\mu'$
  sur $g(l')\subset b\Theta$ s'\'etendent contin\^ument dans $\Theta$ en
  une $(1,0)$-forme holomorphe,  qui ne s'annulle dans aucun
  point.
\item Aucune des fonctions $\overline a+\overline b\overline z_1$,
  $\overline c+\overline d\overline z_1$, $g$ et $z_1$ ne s'annulle sur
  $\overline S_{\nu_0}\cap g^{-1}(\Theta)$.
\end{enumerate}
En consid\'erant le rapport entre  deux mesures $\mu$ et
$\mu'$, on constate que la
restriction de la fonction 
$$h:=g_*\left\{\frac{(\overline c+\overline d\overline
    z_1)(a+bz_1)}{\overline a+\overline b\overline
    z_1}\right\}$$
sur $g(l')$ s'\'etend holomorphiquement sur $\Theta$. Ceci montre que
    la restriction de la fonction 
$$f:=\frac{\overline c+\overline d\overline
    z_1}{\overline a+\overline b\overline
    z_1}$$
sur $\Pi_1(l')$ s'\'etend aussi holomorphiquement dans $\Omega$. Comme la
    matrice ${\cal M}$ est de rang maximal, la fonction $f$ n'est pas
    constante. Par
    cons\'equent, la restriction de $\overline z_1$ s'\'etend
    m\'eromorphiquement dans $\Omega$ car $\overline z_1=(f\overline
    a-\overline c)/(\overline d-f\overline b)$. Un remplacement
    convenable de $l'$, $\Theta$ et $\Omega$ par leurs ouverts permet
    de supposer que cette extension est holomorphe. Soit $\Phi$ une
    application biholomorphe du disque unit\'e $D$ dans
    $\Omega$. On sait que $\Phi$ s'\'etend contin\^ument en une
    application bijective de $\overline D$ dans $\overline\Omega$. 
Alors $\overline \Phi|_{l''}$ s'\'etend holomorphiquement dans
    $D$ o\`u $l'':=\Phi^{-1}(\Pi_1(l'))$. 
En cons\'equence, les restrictions des 
    parties r\'eelle et imaginaire
    $\Re\Phi$ et $\Im\Phi$ de $\Phi$ sur $l''$
    s'\'etendent holomorphiquement dans
    $D$. D'apr\`es le principe de r\'eflexion, les fonctions obtenues se
    prolongent holomorphiquement au voisinage de $l''$. Ceci
    montre que $\Pi_1(l')$ est r\'eelle analytique. En utilisant de
    petits changements des coordonn\'ees, on peut prouver que
    $\Pi_{1,\epsilon}(l')$ contient un ouvert r\'eel analytique pour
    tout 
    $\epsilon$ petit, o\`u
    $\Pi_{1,\epsilon}(z):=z_1+\epsilon_2 z_2+\epsilon_3 z_3$. 
Finalement, $l'$ contient un ouvert r\'eel analytique. 
C'est la contradiction recherch\'ee. 
\end{preuve}
\section{Remarque et questions ouvertes}
La m\^eme m\'ethode de la preuve du th\'eor\`eme 1 permet de
d\'emontrer le r\'esultat suivant:
\begin{theoreme}
Soient $D$ une vari\'et\'e complexe de dimension $p\geq 2$ \`a bord $C^2$
irr\'eductible dans $\mathbb{CP}^n$ et
$f$ une fonction $C^2$ d\'efinie sur $bD$ \`a valeurs dans
$\mathbb{CP}^1$. Soit $V\subset\G(n-p+2,n+1)$ une vari\'et\'e
r\'eelle de codimension $1$ v\'erifiant les conditions suivantes pour
tout $ \nu\in V$
\begin{enumerate}
\item $\bigcup_{\nu\in V}\mathbb{P}^{n-p+1}_\nu$ recouvre un ouvert
  dense de $bD$.
\item $\mathbb{P}^{n-p+1}_\nu$ coupe $bD$ transversalement.
\item La restriction de $f$ sur $bD\cap\mathbb{P}^{n-p+1}_\nu$ se
  prolonge dans $D\cap\mathbb{P}^{1}_\nu$ en une
  fonction m\'eromorphe.
\item Aucun ouvert de $bD\cap\mathbb{P}^{n-p+1}_\nu$ n'est r\'eel analytique.
\end{enumerate}
Alors la fonction $f$ se prolonge m\'eromorphiquement dans $D$.
\end{theoreme}
Ce r\'esultat g\'en\'eralise \cite[th\'eor\`eme
5]{Dinh5}, qui
donne la solution partielle \`a:  
\begin{conjecture} (Globevnik-Stout, \cite{GlobevnikStout}) Soient
  $\Omega\subset\subset D$  deux domaines convexes \`a bord $C^2$ dans
  $\mathbb{C}^n$ et $f$ une
  fonction continue sur le bord de $D$. Supposons que $f$ se prolonge
  contin\^ument en une fonction holomorphe dans $D\cap l$ pour toute
  droite complexe $l$ tangente \`a $b\Omega$. Alors $f$ se prolonge
  holomorphiquement dans $D$.
\end{conjecture}
\begin{conjecture} Soient
  $\Omega\subset\subset D$  des domaines convexes \`a bord $C^2$ dans
  $\mathbb{C}^n$, $K$ un compact convexe de $\mathbb{C}^n$ et $f$ une
  fonction continue sur $bD\setminus K $. Supposons que $f$ se prolonge
  contin\^ument en une fonction holomorphe dans $(D\setminus K)\cap l$
  et aucun ouvert de $(bD\setminus K)\cap l$ n'est r\'eel analytique
  pour toute
  droite complexe $l$ tangente \`a $b\Omega$. Alors $f$ se prolonge
  holomorphiquement dans $D\setminus K$.
\end{conjecture}
\begin{question} Est-ce que  le th\'eor\`eme 3 est valable pour une fonction
  $f$ continue?
\end{question}
\begin{question} Si, dans le th\'eor\`eme 3, $\overline D$ est inclus
  dans $\mathbb{C}^n$ et si la fonction $f$ et ses prolongements dans
  $D\cap\mathbb{C}^{n-p+1}_\nu$ sont tous \`a valeurs dans
  $\mathbb{C}$, 
est-ce que la condition 4 est n\'ecessaire?
\end{question}
\end{document}